%% file: revised-anomalous-relax-FCAA-2015.tex
\documentclass[twoside,reqno,11pt]{fcaa}
\input fcaa_style

\usepackage{hyperref} 
\usepackage{subcaption}
\usepackage{graphicx}

\def\theequation{\arabic{section}.\arabic{equation}}

\def\themycopyrightfootnote{\vspace*{3pt}
\copyright \, Year\,  Diogenes  Co., Sofia
\par  \noindent pp. xxx--xxx
\hfill  \vspace*{-36pt}
  }

\title[Fractional relaxation and oscillations]                                 
{Fractional relaxation \\ and fractional oscillation models \\ involving 
Erd\'elyi-Kober integrals}
\author[M. Concezzi, R. Garra, R. Spigler]{Moreno Concezzi,$^1$ 
Roberto Garra,$^2$ Renato Spigler,$^1$}

\begin{document}

\vbox to 2.5cm { \vfill }


\bigskip \medskip


\begin{abstract}
   We consider fractional relaxation and fractional oscillation equations 
involving Erd\'elyi--Kober integrals. In terms of Riemann--Liouville integrals,
the equations we analyze can be understood as equations with time-varying
coefficients. Replacing Riemann--Liouville integrals with Erd\'elyi--Kober-type
integrals in certain fractional oscillation models, we obtain some more general
integro-differential equations. The corresponding Cauchy-type problems can be 
solved numerically, and, in some cases analytically, in terms of Saigo--Kilbas 
Mittag--Leffler functions. The numerical results are obtained by a treatment 
similar to that developed by K.~Diethelm and N.J.~Ford to solve the
Bagley--Torvik equation. Novel results about the numerical approach to the 
fractional damped oscillator equation with time-varying coefficients are also 
presented.

\medskip

{\it MSC 2010\/}: Primary 34C26; 26A33; 65L05; \, Secondary 326A33; 26A48; 
34A08; 33E12

\smallskip

{\it Key Words and Phrases}: anomalous relaxation; fractional relaxation; 
     fractional oscillations; Erd\'elyi--Kober integrals

\end{abstract}

\maketitle

\vspace*{-16pt}


\section{Introduction}\label{sec:introduction}

\setcounter{section}{1}
\setcounter{equation}{0}\setcounter{theorem}{0}

  In this paper, we consider fractional relaxation and fractional oscillation 
models with time-varying coefficients, involving Erd\'elyi--Kober-type 
integrals. In the recent years, a number of papers have been devoted to the 
analysis of fractional relaxation and oscillations as well as to their 
applications, see, e.g., \cite{gom,ma1} and references therein. For a 
probabilistic approach to fractional relaxation equations, we refer to
\cite{be}. Non-Debye relaxation phenomena in dielectrics are often modeled by 
means of fractional differential equations that generalize the classical 
relaxation equation. In this framework, Mittag--Leffler-type functions play a 
relevant role to describe anomalous relaxation, including special cases such as
the Cole--Cole \cite{Cole}, the Davidson--Cole \cite{Davidson}, and the 
Havriliak--Negami \cite{Havri,Roberto} models (see also \cite{Hilfer} for a 
short review about analytical representations of relaxation functions in 
non-Debye processes). In recent theoretical investigations \cite{ma2,ga}, the 
authors studied fractional relaxation models with time-varying coefficients 
resorting to different approaches. In the first part of this paper, we consider
a further approach in this direction, studying a new fractional model involving 
Erd\'elyi--Kober-type integrals. In this framework, we provide both analytical 
and numerical results.

  In the second part of the paper, we consider a fractional generalization of 
the damped oscillator equation with time-dependent elasticity involving 
Erd\'elyi--Kober-type integrals. Fractional generalizations of the damped 
oscillator equation have been subject of recent research, especially in the 
framework of fractional mechanics; we refer, in particular, to the recent 
books \cite{Kli,Torres} and references therein. Moreover, we observe that a 
connection between fractional diffusion equations involving 
Erd\'elyi--Kober-type operators and the density law of the generalized grey 
Brownian motion has been recently studied in \cite{Pa12}. This fact stresses 
the role of Erd\'elyi--Kober integral operators in the future developements of 
fractional calculus.
 
  As far as we know, only few works concerning analysis as well as numerical 
treatment of fractional damped oscillators with time-dependent coefficients 
can be found in the literature. We will show that some generalizations of the 
classical equations of mechanics by means of Erd\'elyi--Kober-type integrals 
imply the analysis of integro-differential equations with time-varying 
coefficients involving Caputo-type fractional derivatives. This topic seems to 
be promising and interesting both, for numerical and analytical studies of 
fractional differential equations with variable coefficients, and their 
applications in mechanics. 

%
\section{Fractional relaxation models involving Erd\'elyi--Kober-type integrals}
\label{sec:1}

\setcounter{section}{1}
\setcounter{equation}{0}\setcounter{theorem}{0}

  A reasonable generalization of the classical relaxation model can be obtained
replacing the ordinary time derivative with the the Caputo fractional 
derivative, defined by (see \cite{ma1})
\begin{equation} \label{Capu}
      D^{\alpha}_{t}f(t) = \frac{1}{\Gamma(m-\alpha)}
                 \int_0^{t}(t-s)^{m-1-\alpha} \frac{d^m}{ds^m}f(s) \, ds,
\end{equation}
where $\alpha > 0$, $m = \lceil \alpha \rceil$, and $f \in AC^m[a,b]$. 
Therefore, the so-called fractional relaxation equation is simply given by
\begin{equation}\label{1.1}
       D^{1-\alpha}_{t}u(t)
     = -\lambda u(t), \quad \lambda > 0, \ t \geq 0, \ \alpha \in [0,1).
\end{equation}
In this case, the fractional model describes an anomalous relaxation process 
with a Mittag--Leffler decay, that is with an asymptotically power law behavior.

  Since for absolutely continuous functions, we have
\begin{equation}
       J_t^{1-\alpha}  D^{1-\alpha}_t u(t) = u(t)-u(0),
\end{equation}
where 
\begin{equation}
   J_t^{\beta} f(t) = \frac{1}{\Gamma(\beta)} \int_0^t(t-s)^{\beta-1} f(s) \, ds,
\end{equation}
is the Riemann-Liouville fractional integral of order $\beta > 0$ (see, 
\cite{e-p}, e.g.), we can write \eqref{1.1}, in the equivalent integral form,
\begin{equation}\label{1.2}
           u(t)-u(0) = -\lambda J_t^{1-\alpha} u(t). 
\end{equation}

  In this section we consider a generalization of the fractional relaxation 
equation in the integral form \eqref{1.2}, where the Riemann--Liouville 
integral is replaced by the Erd\'elyi--Kober-type fractional integral (see for 
example \cite{mc2, sn1, sn2})
\begin{equation}\label{ek}
     I_m^{\eta, \alpha} f =
   \frac{t^{-m\eta-m\alpha}}{\Gamma(\alpha)} \int_0^t(t^m-s^m)^{\alpha-1} s^{m\eta} 
       f(s) \, d(s^m), \quad \eta \geq 0, \ m > 0.
\end{equation}

  By this, we will be able to obtain some general results concerning anomalous
relaxation with time-dependent coefficients and hereditary effects, including 
both, classical (ordinary) and fractional relaxation models as special cases.
 
  More general Erd\'elyi--Kober operators, with $\beta > 0$ instead of $m$ in 
\eqref{ek}, have been studied, and the corresponding ``Erd\'elyi--Kober 
fractional derivatives'' were introduced in \cite{kiry}, and more recently 
analyzed, e.g., in \cite{ltr}. Applications of such operators in the framework 
of fractional mechanics should be the subject of further research.

  For simplicity, in \eqref{ek} we will assume $m = 1$, being $\alpha \in 
(0,1)$. In fact, in this section we consider purely relaxation phenomena. 
Moreover, we assume that the relaxation coefficient there follows a power-law 
behavior, i.e., $\lambda(t) = \lambda t^{1-2\alpha}$. We are interested to 
understand the interplay between dynamical coefficients and 
Erd\'elyi--Kober-type integrals in generalized mechanical models.
  Therefore, our first model is based on the following generalization of 
equation \eqref{1.2}
\begin{align} 
   \nonumber 
  u(t) - u_0 & = -\frac{\lambda}{t^{2\alpha-1}}I_1^{\eta,1-\alpha}u(t)  \\
  \nonumber  & =-\frac{\lambda \, t^{-\alpha-\eta}}{\Gamma(1-\alpha)} 
                  \int_0^t (t-s)^{-\alpha}s^{\eta} u(s) \, ds  \\
             & = -\lambda \, t^{\alpha-\eta}J_t^{1-\alpha} \left(t^{\eta}u(t) \right).
                                                          \label{cia}
\end{align}
We observe that, for $\eta = 0$ and $\alpha = 0$ we recover the classical 
relaxation model. Hence, the effect of introducing Erd\'elyi--Kober-type 
integrals in the fractional relaxation equation is essentially that of providing
power law time-varying coefficients. We will show that this integral equation 
can be treated in a similar way as that considered in \cite{ma2}.

  Equation \eqref{cia} can be handled in two equivalent ways. The first one is 
based on the solution of the following system of coupled integro-differential 
equations
\begin{equation}\label{capu}
  \begin{cases}
      u(t) - u_0 = -\lambda t^{-\eta-\alpha} D^{\alpha}_t g(t) \\
      D_t g = t^{\eta} u(t),
  \end{cases}
\end{equation}
which is clearly equivalent to \eqref{cia} (here $D_t g \equiv D^1_ t := dg/dt$
and $g(0) = 0$). This approach is interesting in view of a numerical 
integration of equation \eqref{cia}, and is suggested by a method developed by 
Diethelm and Ford in \cite{Dieth02,Dieth04} to solve numerically the well-known 
Bagley--Torvik equation \cite{bagley}.

  The second approach, based on simple analytical manipulations, allows to 
reduce equation \eqref{cia} to that considered recently in \cite{ma2}.
  Indeed, applying the first order derivative to both sides of \eqref{cia}, we 
obtain
\begin{equation}
  \frac{d}{dt} \left[u(t) + \lambda t^{-\eta-\alpha}J_t^{1-\alpha}(t^{\eta} u) 
      \right] = 0.
\end{equation} 
We will study the case  
\begin{equation}\label{cici}
      t^{\eta+\alpha}u(t) + \lambda J_t^{1-\alpha}(t^{\eta}u) = 0.
\end{equation}
Let us recall that
\begin{equation}
       D^{1-\alpha}_t J_t^{1-\alpha} f(t) = f(t),
\end{equation}
and define $r(t) := t^{\eta+\alpha} u(t)$. Then, taking the Caputo derivative of 
order $1 - \alpha$ of both sides of \eqref{cici}, we obtain the equation
\begin{equation} \label{cidi}
       D^{1-\alpha}_t r(t) = -\lambda t^{-\alpha} r(t).
\end{equation}
Such equation has recently been studied by Capelas de Oliveira et al. in 
\cite{ma2}, and its explicit solution can be written in terms of the so-called 
Saigo--Kilbas Mittag--Leffler function, for $\alpha \in (0, 1/2)$ (see the 
recent monograph \cite{Rogosin} for more details). In particular, here we 
recall the following useful result (see, e.g., \cite{e-p}, pp.~232--233).
\begin{lemma}  
The function
\begin{align}\label{saigo}
           E_{\alpha,1+\frac{\beta}{\alpha},\frac{\beta}{\alpha}} & \left(-
                \lambda t^{\alpha + \beta} \right)
               = 1 + \sum_{k=1}^\infty (-\lambda)^k t^{k \left( \alpha + \beta \right)}
    \prod_{j=0}^{k-1} \frac{\Gamma \left( \alpha \left( j+j \frac{\beta}{\alpha}
   + \frac{\beta}{\alpha} \right)+1 \right)}{\Gamma \left(
    \alpha \left( j + j \frac{\beta}{\alpha} + \frac{\beta}{\alpha} 
           + 1 \right) +1 \right)},       
\end{align}
solves the following fractional Cauchy problem:
\begin{equation} \label{odefr}
  \begin{cases}
     D^{\alpha}_t y \left(t \right) = -\lambda t^\beta y \left( t \right), &
         t \ge 0, \: \alpha \in \left( 0, 1 \right], \: -\alpha < \beta 
          \leq 1-\alpha, \\
     y \left( 0 \right) = 1.
   \end{cases}
\end{equation}
\end{lemma}

  We recall that, according to Kilbas et al. \cite{e-p} (Remark~4.8, p.~233), 
uniqueness of the solution to the Cauchy problem in \eqref{odefr} has been 
established only for $\beta \geq 0$.

  In view of Lemma 1.1, a solution to equation \eqref{cidi}, for $\alpha \in 
(0,1/2)$, is given by 
\begin{equation}
   r(t) = E_{1-\alpha, 1-\frac{\alpha}{1-\alpha}, -\frac{\alpha}{1-\alpha}}
            \left(-\lambda t^{1-2\alpha}\right),
\end{equation}
and going back to the original function, we have that
\begin{equation}\label{sol0}
     u(t) = \frac{ E_{1-\alpha, 1-\frac{\alpha}{1-\alpha}, -\frac{\alpha}{1-\alpha}}
            \left(-\lambda t^{1-2\alpha}\right)}{t^{\alpha+\eta}}.
\end{equation}
We observe that, for $\alpha = \eta = 0$, equation \eqref{sol0} becomes
\begin{equation}
     u(t) = E_{1, 1, 0}\left(-\lambda t\right) = e^{-\lambda t},
\end{equation}
which coincides with the solution of the classical relaxation equation with 
initial condition $u(0) = 1$.

  The study of the conditions for the complete monotonicity of a solution to 
the fractional Cauchy problem in \eqref{odefr} is a central topic for the 
physical sense of the anomalous relaxation models, since it ensures that in 
isolated systems the energy decays monotonically (see \cite{Hanyga} for a 
complete discussion about the \emph{physical realizability} of the
\emph{viscoelastic models}). 

  In \cite{ma2}, the authors have discussed the complete monotonicity of the 
so-called Saigo--Kilbas Mittag--Leffler function \eqref{saigo}. In particular, 
they proved that the solution \eqref{saigo} is completely monotonic  provided 
that $\alpha \in (0,1]$ and $\beta \in (-\alpha, 1 - \alpha]$. In our case,
these conditions correspond to the assumption that $\alpha \in (0,1/2)$, that 
ensures the complete monotonicity of the solution. We refer to \cite{rgm} and 
references therein for a recent discussion concerning the relevance of 
completely monotonic functions in dielectrics.

%
\subsection{Numerical results}

\begin{figure}[htbp]
\centering
\includegraphics[scale=0.3]{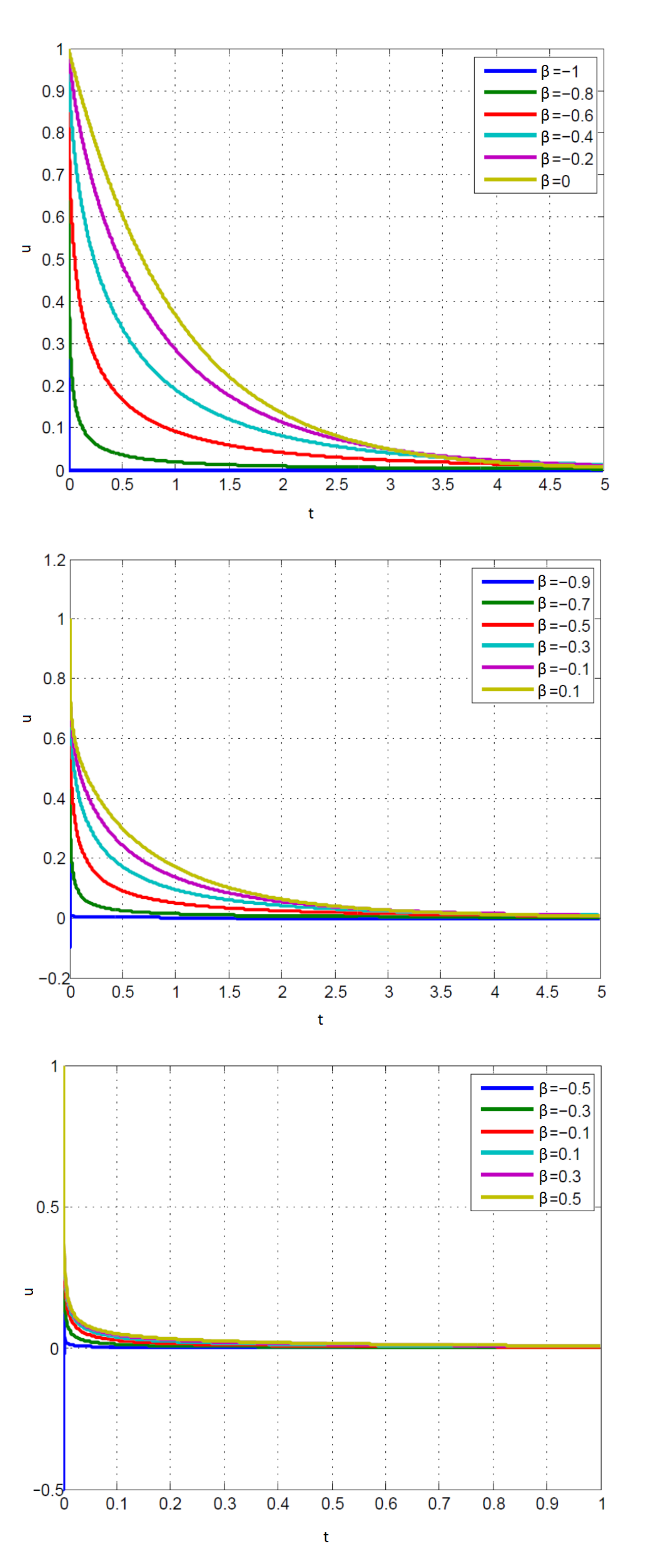}
\caption{Numerical solution of (\ref{odefr}) for $\alpha = 0.9$, $\lambda = 1$,
initial discretization step $h = 0.001$ and several values of $\beta$, computed 
for $t \in [0,5]$.}\label{fig_3_4_5} 
\end{figure}

  In this section, we present some numerical results concerning the solution of
the fractional Cauchy problem in \eqref{odefr}. In Fig.~\ref{fig_3_4_5}, we
show the numerical solutions to such problem, for different values of $\alpha 
\in (0,1]$, and varying the parameter $\beta$, with $- \alpha < \beta \leq 1 
- \alpha$

  The decay of the response function, for a fixed $\alpha$, is faster than 
exponential for short times, and power-law-like asymptotically, as pointed out 
also in \cite{ma2}. The decay is clearly faster for lower values of $\alpha$. 

\begin{figure}[h!]
\centering
\includegraphics[width=10cm]{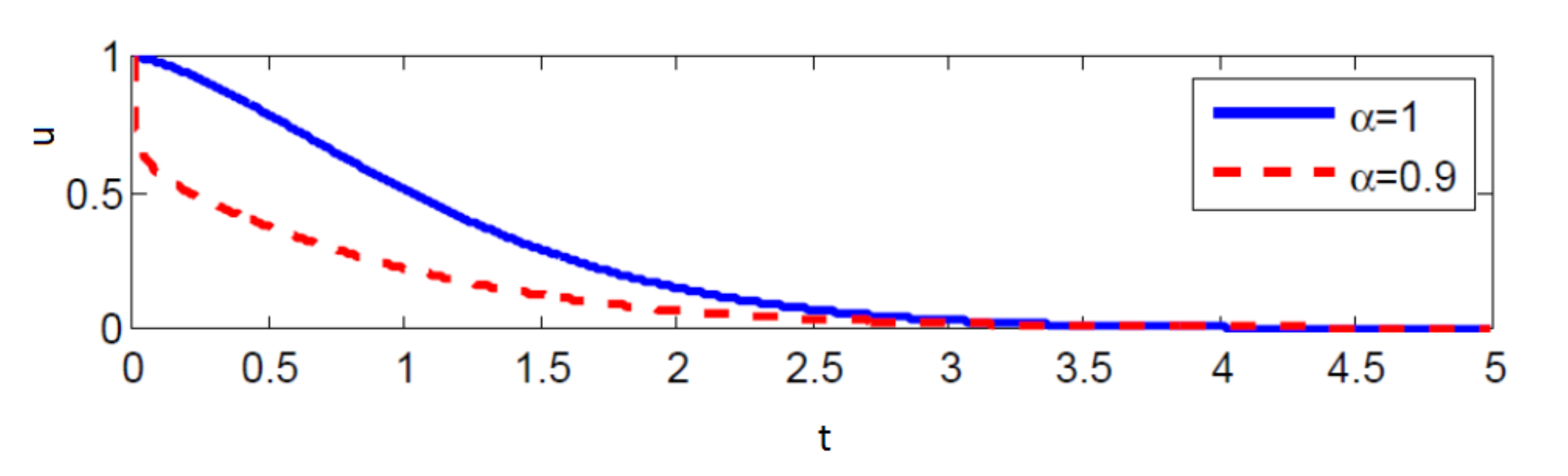}
\caption{Numerical solution of (\ref{odefr}) for $\alpha = 1$ and $\alpha = 
0.9$, $\lambda = 1, \beta = 0.5$, and initial discretization step $h = 0.001$, 
computed for $t \in [0,5]$.}
\label{fig_2}
\end{figure}

  In Fig.\ref{fig_2}, we observe that, for a fixed value of $\alpha$ ($\alpha = 
0.9$) for short times the solution to equation \eqref{saigo} decays faster with 
respect to the classical exponential relaxation decay.

  On the other hand, for a fixed value of $\alpha$, the role played by the 
parameter $\beta$ is clear: increasing $\beta$ the solution decays faster (see 
Fig.\ref{fig_3_4_5}). In order to understand the physical meaning of the 
fractional relaxation model described by the Cauchy problem in \eqref{odefr}, 
we observe that, for $\alpha = 1$, it reduces to
\begin{equation} \label{oded}
  \begin{cases}
      D_t y \left(t \right) = -\lambda \, t^\beta y \left( t \right), &
          t \ge 0,  \,  \beta > -1 \\
       y\left( 0 \right) = 1,
   \end{cases}
\end{equation}
that is to a relaxation equation with the time-dependent coefficient 
$\lambda(t) := \lambda \, t^{\beta}$, whose solution is given by
\begin{equation}\label{sold} 
      \nonumber 
    y(t) = \exp{\left\{- \lambda \, \frac{t^{\beta+1}}{\beta+1} \right\} }.
\end{equation}
  From this special case, we can understand the role of the parameter $\beta$ 
in the solution, as confirmed by the numerical results.

  In Fig.s~\ref{fig_2} and ~\ref{fig_3_4_5}, we plotted the numerical solution 
of (\ref{odefr}) for several values of the parameters. The numerical method we 
used is an adaptive improvement of the predictor-corrector method earlier 
introduced in \cite{Dieth99}, and developed in \cite{Conc}. 
  Starting from an arbitrary discretization step, $h$, we chose locally a step 
size inversely proportional to the size of the (classical) derivative of the 
solution being computed, so that, at the $i$-th step, the time step 
$$
      h_i := \frac{c_i \, h}{|u(t_i) - u(t_{i-1})|},
$$
$i = 1, 2, \ldots, N$, $t_i = \sum_{j=1}^i h_j$, $t_0 = 0$, is used. Here $N$ is
the number of the integration nodes, and
$$
       c_i = \frac{u^\prime(t_{i-2}) - u^\prime(t_{i-3})}
                {u^\prime_\alpha(t_{i-1}) - u^\prime(t_{i-2})},
$$
where $u^\prime$ denotes the (classical) first derivative of $u^\prime$, see 
\cite{Conc}.


\begin{figure}[h!]
\begin{center}
\includegraphics[width=10cm]{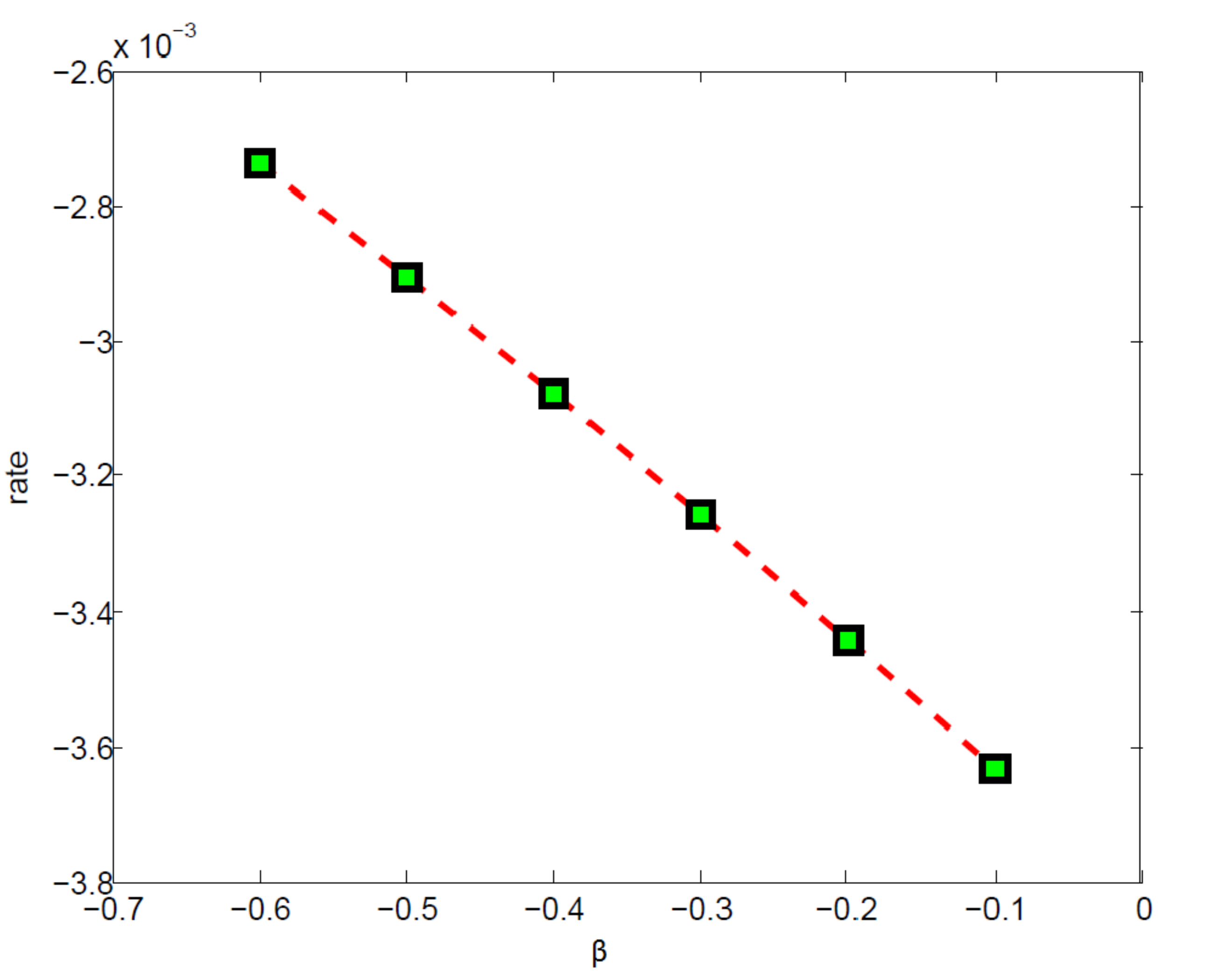}
\end{center}
\caption{Rate of convergence, represented by $\log|u^{\prime}/u|$, for the 
numerical solution of (\ref{odefr}) for $\alpha = 0.9$, $\lambda = 1, \beta 
\in [-0.6,-0.1]$, initial discretization step $h = 0.1$, computed fro $t \in 
[0,500]$, in a $log-log$ scale.}
\label{rate_1}
\end{figure}

 In Fig.~\ref{rate_1}, we show the convergence rate for the numerical solution 
of (\ref{odefr}), plotting $\log|u^{\prime}/u|$ vs $\beta$. The graph exhibits a 
linear {\em negative} relation between the values of $\beta$ and the rate of 
convergence defined as before.

%

%
\section{Damped fractional oscillator involving Erd\'elyi--Kober-type integrals}
\label{sec:2}

\setcounter{section}{2}
\setcounter{equation}{0}\setcounter{theorem}{0}

  In this section, we consider a fractional generalization of the damped 
oscillator equation, involving Erd\'elyi--Kober-type integrals. There are many 
investigations nowadays about fractional damped oscillators involving Caputo 
derivatives, in the framework of the fractional mechanics (see for example 
\cite{duang,Kli,Sca} and references therein). Here we consider a special form
of the fractional damped equation with time-varying coefficients, containing 
Erd\'elyi--Kober-type integrals, and we assume the elastic term to be time 
dependent. Hereafter, $D_t$ will denote the first order ordinary time 
derivative. A Caputo-type fractional generalization of the damped oscillator 
equation is thus given by
\begin{equation}\label{h1}
         D^{\alpha}_tD_t x(t) + \lambda D_t x(t) + k x(t) = 0,
\end{equation} 
where $\lambda > 0$ represents the viscosity coefficient and $k > 0$ is the 
elasticity constant. Since ${}^C D^{\alpha}_t x(t) = J_t^{1-\alpha} D_t x(t)$, 
\eqref{h1} can also be written as
\begin{equation}\label{h2}
     J_t^{1-\alpha}\left(D^2_t x\right)(t) + \lambda D_t x(t) + k x(t) = 0. 
\end{equation}
Our purpose is to consider a further, more general model of damped oscillator
with memory effects and time-varying coefficients, replacing in \eqref{h2} the 
Riemann--Liouville integral with the Erd\'elyi-Kober-type integral 
\begin{equation} \label{sec}
    \left( I_1^{\eta,1-\alpha} u \right)(t)
  = \frac{t^{\alpha-1-\eta}}{\Gamma(1 - \alpha)} \int_0^t(t-s)^{-\alpha} s^{\eta}
       u(s) \, ds.
\end{equation}
Moreover, we will assume that the elasticity term depends on time, and, in 
particular, $k(t)= k_0/t$, where we set $k_0 := \lambda(\eta + 1 - \alpha) > 
0$. 
We thus obtain the following fractional damped oscillator equation with 
time-varying coefficients,
\begin{equation}\label{h3}
       \left[D^{\alpha}_t t^{\eta}D_t + \lambda t^{\eta+1-\alpha} 
       D_t +  k_0 t^{\eta-\alpha} \right] u(t) = 0.
\end{equation}
Note that, for $\eta = 0$ and $\alpha = 1$, \eqref{h3} reduces to the classical
damped oscillator equation with time-dependent elastic term. At best of our 
knowledge no papers can be found in the literature concerning fractional damped
oscillator equations with time-varying coefficients. Looking at \eqref{h3}, we 
can see that, replacing the Riemann-Liouville integral with the Erd\'elyi-Kober
integral, power law time-variable viscosity and elasticity coefficients appear.

  The physical meaning of having such a combination (in $k_0$) of the original 
constants, $\alpha$ and $\eta$, is that elasticity, damping, and viscosity 
effects turn out to be coupled, that is, we may expect that damping ratio and 
(local) oscillation frequency (if any) are linked. However, we will see from 
the numerical simulations that oscillatory effects strongly depend on the order
$\alpha$ of the fractional derivatives involved in the model equation. 

  We now show how this choice is useful from the mathematical point of view.
Equation \eqref{h3} can be rewritten, decoupling it into the two auxiliary 
coupled equations
\begin{equation}\label{h4}
  \begin{cases}
          D^{\alpha}_t g(t) = -\lambda t^{\eta+1-\alpha}u(t),  \\
         t^{\eta}D_t u(t) = D_t g(t).
  \end{cases}
\end{equation}  
This can be shown just by taking the fractional derivative of order $\alpha$ of
both sides of the second equation of the system \eqref{h4}, and replacing there
$D_t^{\alpha} g$ as given by the first equation, recalling that $k_0 = \lambda 
(\eta + 1 - \alpha)$. We obtain
\begin{equation}
       D^{\alpha}_t D_t g(t) = D_t  D^{\alpha}_t g(t)
         = -\lambda D_t \left(t^{\eta+1-\alpha} u(t) \right)
         =  D^{\alpha}_t\left(t^{\eta} D_t u(t) \right),
\end{equation}
where we exploited the fact that, being $g(0) = 0$, ordinary and fractional 
derivatives commute. Moreover, in this case the Riemann--Liouville and the 
Caputo fractional derivatives of $g(t)$ coincide (see, e.g., 
\cite[pp.~73-74]{pod}. This seems to be an interesting result, since equation 
\eqref{h3} is a fractional differential equation with time-varying 
coefficients, classically related to fractionally damped oscillation.    

  The previous equation can be written as a ``hybrid'' system, in the sense 
that it involves both, classical and fractional derivatives. Here, the range of
the parameters is $0 < \alpha \leq 1$, $\lambda > 0$, and $\eta \geq 0$. 
Writing $v(t) := D_t u(t)$, we obtain
\begin{equation}
  \left\{
      \begin{array}{ll}
     D_t u(t) = v(t)  \\
     D_t^{\alpha} \left( t^{\eta} v \right)(t) 
              = - \lambda t^{\delta} v - \mu t^{\gamma} u,
       \end{array}   \right\},
\end{equation}
where we set, for short, $\delta := \eta + 1 - \alpha$, $\mu := \lambda (\eta 
+ 1 - \alpha)$, and $\gamma := \eta - \alpha$, and assume as initial conditions
$u(0) = v(0) = 1$. Setting $w(t) := t^{\eta} v(t)$, we have 
\begin{equation}\label{hybrid_system}
  \left\{
       \begin{array}{ll}
     D_t u(t) = t^{-\eta} w(t)  \\
     D_t^{\alpha} w(t) = - \lambda t^{\delta - \eta} w(t) - \mu t^{\gamma} u(t),
       \end{array}   \right.
\end{equation}

  This approach, adopted to solve numerically \eqref{h3}, is suggested again by
the method of Diethelm and Ford in \cite{Dieth02}, extended according to the
adaptive pattern developed in \cite{Conc}), to solve the celebrated 
Bagley--Torvik equation \cite{bagley}, 
\begin{equation} \label{BT}
         A D_t^2 y + B {}^C D_t^{3/2} y + C y = f(t),
\end{equation}
where $A \neq 0$, $B$, and $C$ are real constants. We recall that in 
\cite{Dieth02} equation \eqref{BT} was rewritten as
\begin{equation} \label{BT_system}
  \left\{
      \begin{array}{llll}
           D_t^{1/2} y_1 = y_2  \\
           D_t^{1/2} y_2 = y_3  \\
           D_t^{1/2} y_3 = y_4  \\
           D_t^{1/2} y_4 = A^{-1} [-C y_1 - B y_4 + f(t)],  
       \end{array} \right.
\end{equation} 
with the initial conditions $y_1(0) = y_0$, $y_2(0) = 0$, $y_3 = y_0'$, $y_4(0)
= 0$. However, using Gr\"{u}wald-Letnikov discretization, we face the problem 
of having unevenly spaced nodes, which problem can be solved resorting to some 
interpolation. The latter can be merely linear, since the additional error 
introduced is of the second order, while the order of the method of 
\cite{Dieth02} is of order $3/2$. Quadratic interpolation was also tried, 
which introduces third-order (hence negligible) errors, at the price of some 
little complications.

%
\section{Numerical results} \label{sec:3}

\setcounter{section}{3}
\setcounter{equation}{0}\setcounter{theorem}{0} 

  In this section, we present  some results for the numerical solution of 
\eqref{hybrid_system}. We first discuss the role played by the parameters 
$\alpha$, $\lambda$, and $\eta$ entering equation \eqref{hybrid_system}. The 
fractional order of the derivative has an intrinsic damping effect on 
oscillatory phenomena, as discussed for example in \cite{PhysicaA}. In order to
understand the role of the parameters $\lambda$ and $\eta$, we observe that for
$\alpha = 1$, \eqref{hybrid_system} becomes
\begin{equation} \label{og}
    \left[D_t^2 + \left(\frac{\eta}{t} + \lambda\right)D_t
         + \frac{\lambda\eta}{t} \right]u(t) = 0,
\end{equation}
where we recall that we assumed $\lambda > 0$ and $\eta \geq 0$. This is a 
damped oscillator with time-dependent viscosity and elasticity terms. We remark
that in our fractional model, time-variable viscosity, damping, and elasticity 
terms are coupled. This implies that damping and (instantaneous) oscillation 
frequency are linked. When $\eta = 0$ and $\alpha = 1$, \eqref{hybrid_system} 
reduces to the classical relaxation equation with constant coefficients.

  Equation \eqref{og} can actually be solved explicitly, see, e.g., 
\cite{Kamke,Polyanin,Tricomi}. A special solution is given by $u(t) = 
e^{- \lambda t}$, as one can immediately check, but the general solution is also 
available:
$$
  u(t) = e^{- \lambda t} \left\{c_1 + c_2 \int^t v^{- \eta} e^{\lambda v} dv \right\},
$$
$c_1$ and $c_2$ being two arbitrary constants. If $c_2 \neq 0$,we can see that
$$
  u(t) \sim c_2 \frac{1}{\lambda} \, t^{- \eta},  \ \mbox{as} \  t \to +\infty,
$$   
for $\eta > 0$, being $\lambda > 0$, that is, $u(t)$ decays monotonically as 
$t$ gets large. The oscillatory behavior of the general solution can be seen
proceeding as follows. Setting $u := p(t) v$ with $p(t) = t^{- \eta/2} 
e^{- \lambda t/2}$, $v$ will satisfy the equation
\begin{equation} \label{eq-for-v}
      v'' + q(t) v = 0,
\end{equation}
where
\begin{equation} \label{def-of-q}
      q(t) := 
  \frac{\eta(2 - \eta)}{4 t^2} + \frac{\lambda \eta}{2 t} - \frac{\lambda^2}{4}.
\end{equation}
Now, as is well known, where $q(t) \geq k^2 > 0$, for any arbitrary but fixed 
constant $k$, the solution is oscillatory in the interval, $(t_-(k), t_+(k))$,
where
$$
     t_{\pm} = \frac{\lambda \eta \pm \sqrt{2 \eta \lambda^2 + 4 k^2 
              \eta (2 - \eta)}}{\lambda^2 + 4 k^2}.
$$

For small $k$'s we have 
$$
     t_+(k) - t_-(k) \approx 2 \frac{\sqrt{2 \eta}}{\lambda} 
         \left[1 - \frac{k^2}{\lambda^2}(2 + \eta) \right],
$$
and hence about $t_+(k) - t_-(k) \approx 2 \sqrt{2 \eta}/\lambda$. By the 
Sturm's comparison theorem, the distance between consecutive zeros of $u(t)$ is 
less than $\pi/k$. In fact, by such a theorem, between any two consecutive 
zeros of any nontrivial solution to $u''+ k^2 u = 0$, there exists at least one 
zero of the solution to (\ref{eq-for-v})-(\ref{def-of-q}). Since $t_+(k) 
- t_-(k)$ (for small $k$'s) decreases roughly according to 
$\sqrt{\eta}/\lambda$, we infer that several oscillations occur in the interval
$(t_-, t_+)$.

  On the other hand, the numerical results below show the ``physical'' role of 
the fractional order of derivatives, $\alpha$, which is to introduce a rather 
strong oscillatory behavior in the solutions, for some special combinations of 
the parameters $\lambda$ and $\eta$.

\begin{figure}[h!]
\begin{center}
\includegraphics[scale=0.30]{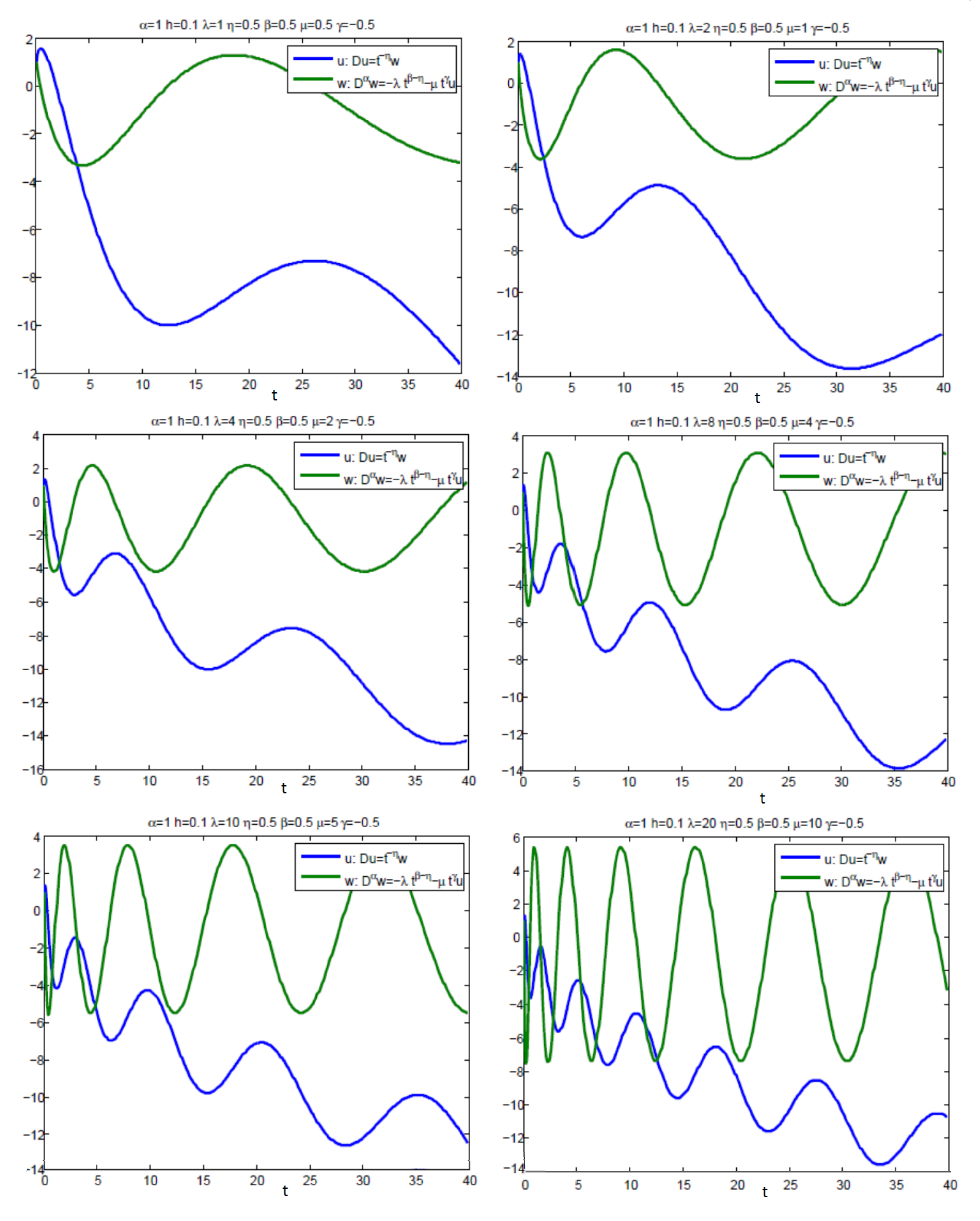}
\end{center}
\caption{Numerical solution of (\ref{hybrid_system}) for $\alpha = 1.0, \eta =
0.5, \mu = 0.5, \beta = 0.5, \gamma = -0.5$, and different values of $\lambda$,
initial discretization step $h = 0.1$, computed for $t \in [0,40]$.}
\label{lambda_more}
\end{figure}

  Generally speaking, varying the parameters $\lambda$, $\eta$, and $\alpha$, 
in analogy with the well-known solutions of the classical damped oscillator 
equation, the solutions can be either {\em overdamped} (i.e., they may have a 
purely monotonic decay), or {\em underdamped} (i.e., they may exhibit an 
oscillatory decay). In our numerical solutions, we observe both behaviors, 
according to the physical parameters. In Fig.~\ref{lambda_more}, we plotted the
numerical solution of (\ref{hybrid_system}), for several values of $\lambda$, 
keeping fixed the other coefficients. We observed an underdamped-like 
oscillatory decay. In Fig.~\ref{more_alpha_fixed_eta_h_001}, we consider the 
solutions for several values of $\alpha $ and $\lambda$, for a fixed value of 
$\eta$. In this case, for a fixed $\alpha$, changing $\lambda$, the behavior of
the solution switches from an overdamped-like monotonic decay to an oscillatory
decay pattern, for increasing values of $\lambda$, as expected. A similar less 
pronounced switch between different behaviors can be observed in 
Fig.~\ref{more_alpha_fixed_lambda_h_001}, varying $\eta$ and $\alpha$.

\begin{figure}[h!]
\begin{center}
\includegraphics[width=7.5cm]{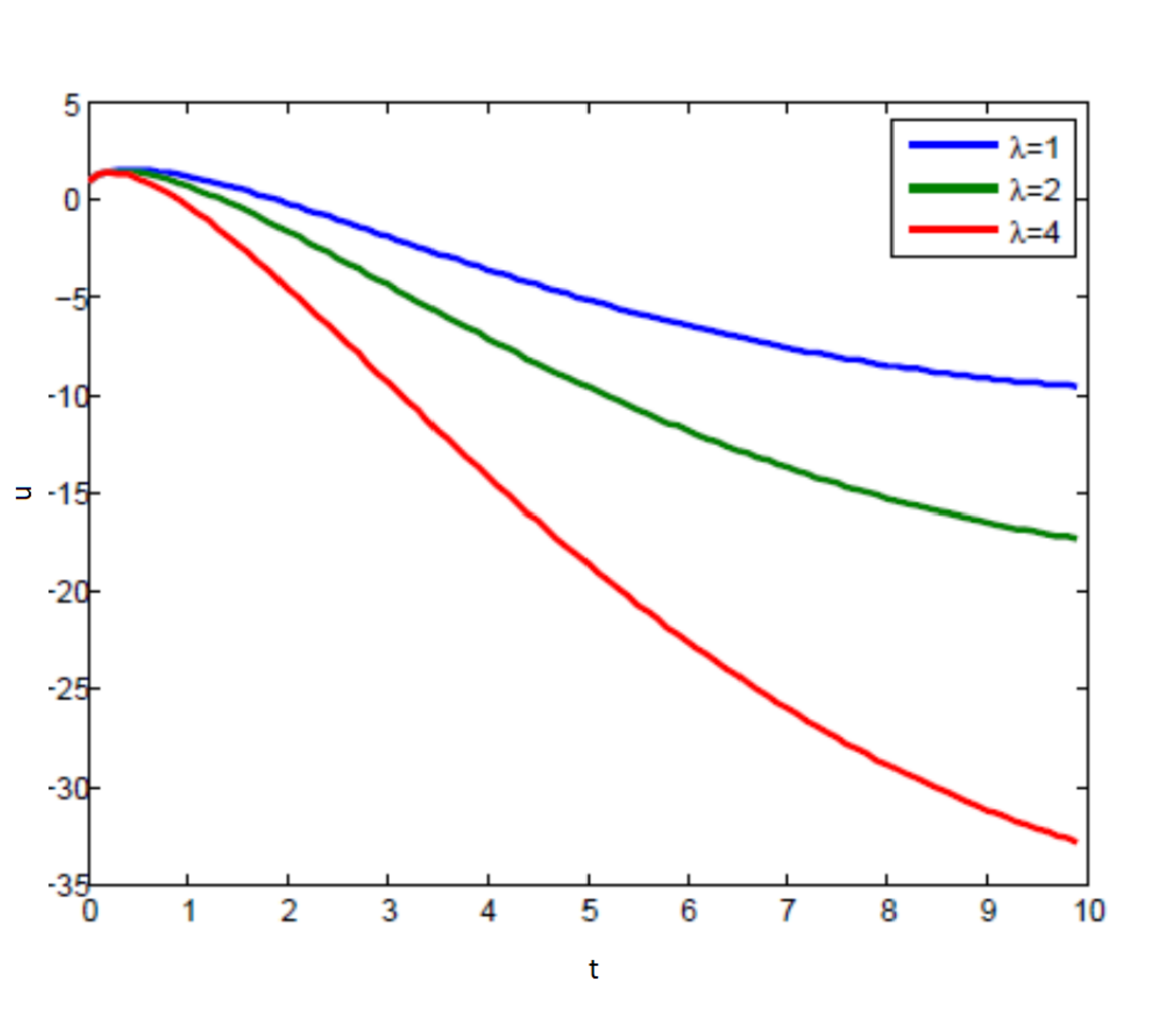}
\end{center}
\caption{Numerical solution of (\ref{hybrid_system}) for $\alpha = 1.0$, $\eta
= 0.5, \mu = 0.5, \beta = 0.5, \gamma = -0.5$, and different values of 
$\lambda$, initial discretization step $h = 0.1$, computed for $t \in [0,10]$.}
\label{more_lambda}
\end{figure}

\begin{figure}[h!]
\begin{center}
\includegraphics[width=8cm]{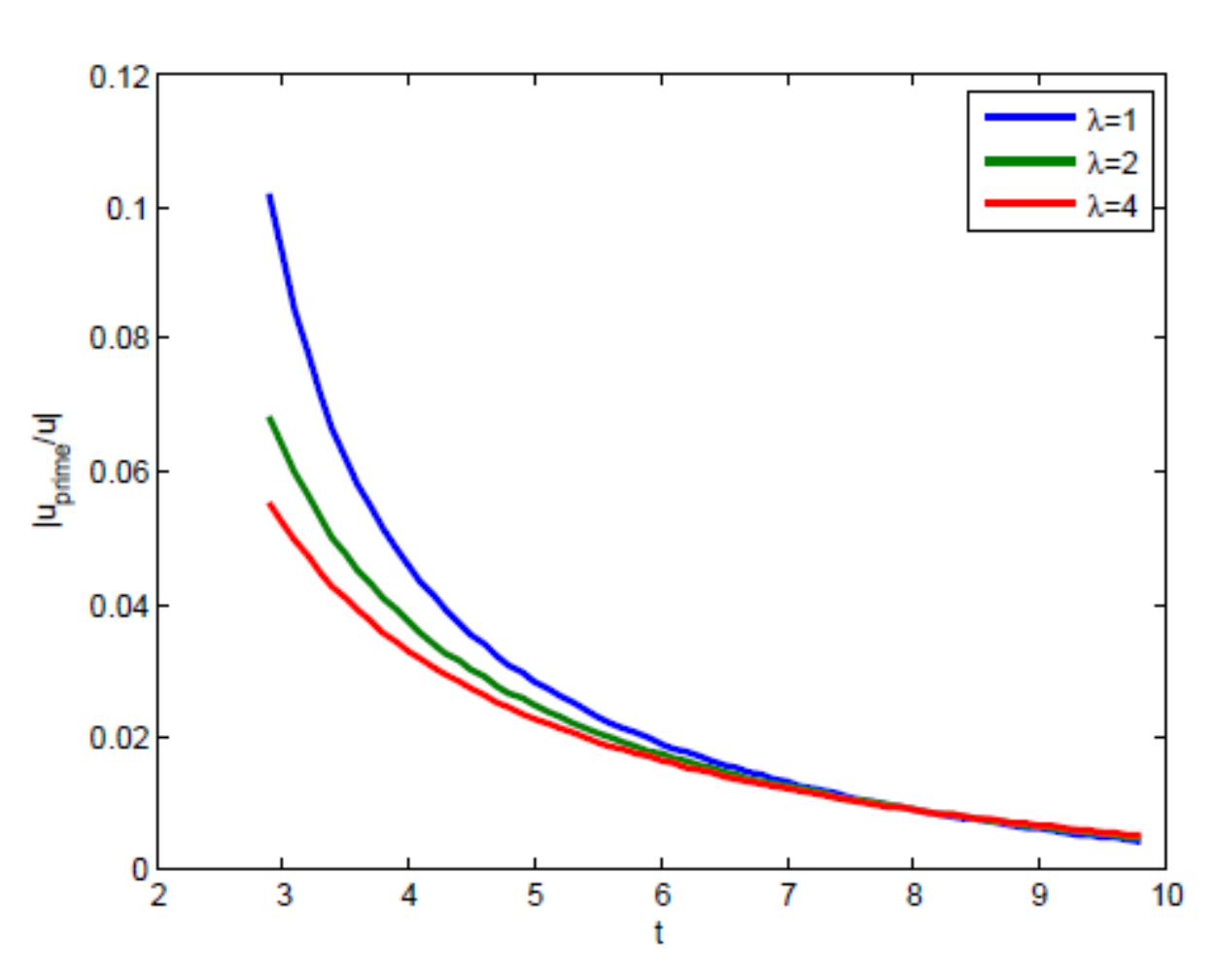}
\end{center}
\caption{Numerical solution of (\ref{hybrid_system}) plotted as $|u^\prime/u|$ 
for $\alpha = 1.0$, $\eta = 0.5, \mu = 0.5, \beta = 0.5, \gamma = -0.5$, and 
several values of $\lambda$, initial discretization step $h = 0.1$, computed 
for $t \in [2,10]$.}
\label{rate}
\end{figure}

  In Fig.~\ref{rate}, we plotted $|u^\prime/u| \ \mbox{vs} \ t$, for $\lambda =
1, 2, 3$, correspondingly to the numerical solution of (\ref{hybrid_system}).

\begin{figure}[h!]
\begin{center}
\includegraphics[width=10cm]{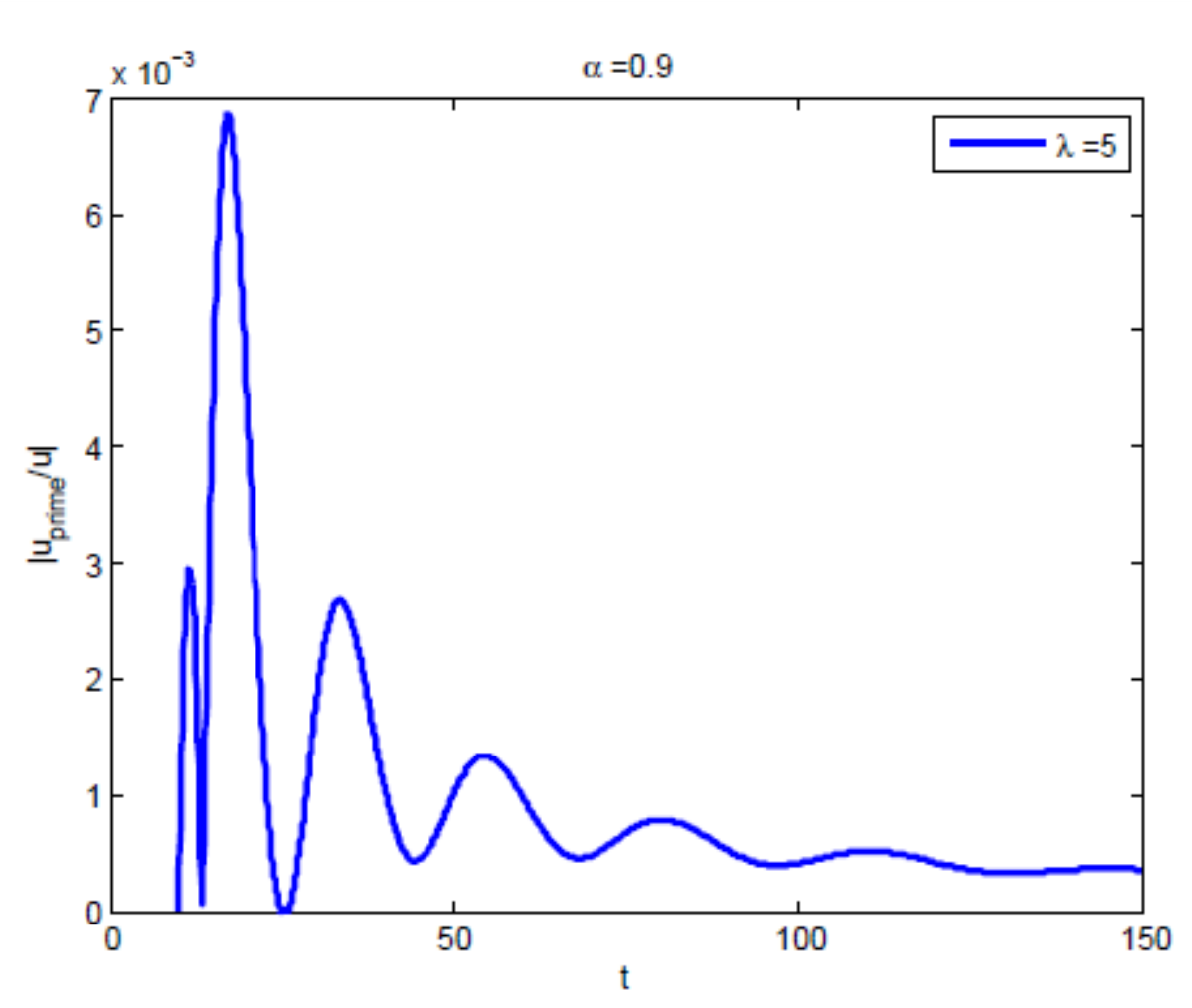}
\end{center}
\caption{Numerical solution of (\ref{hybrid_system}), plotted as $|u^\prime/u|$,
for $\alpha = 0.9$, $\lambda = 5, \eta = 0.5, \mu = 0.6,\beta = 0.6, \gamma =
-0.4$, the initial discretization step $h = 0.1$, computed for $t \in [0,150]$.}
\label{abs_u_prime__over_u_vs_t}
\end{figure}

  In Fig.~\ref{abs_u_prime__over_u_vs_t}, we plotted $|u'/u| \ \mbox{vs} \ t$, 
for $\lambda = 5$ and large $t$, correspondingly to the numerical solution of 
(\ref{hybrid_system}).

\begin{figure}[h!]
\begin{center}
\includegraphics[scale=0.3]{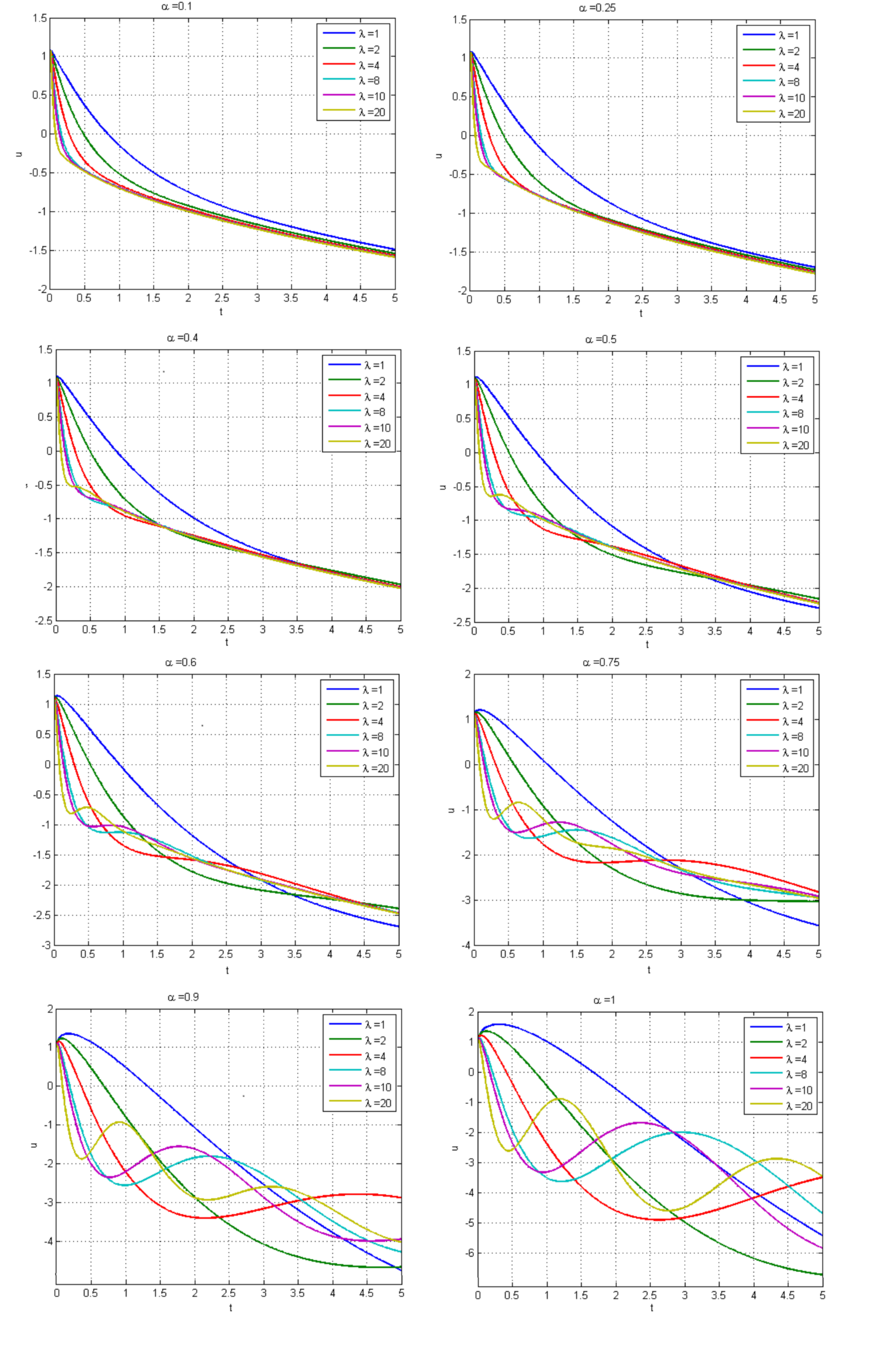}
\end{center}
\caption{Numerical solution of (\ref{hybrid_system}), plotted for several 
values of $\alpha \in (0,1]$ and $\lambda$, with fixed $\eta = 0.5$, 
initial discretization step $h = 0.001$, computed for $t \in [0,5]$.}
\end{figure} \label{more_alpha_fixed_eta_h_001}
\begin{figure}[h!]
\begin{center}
\includegraphics[scale=0.43]{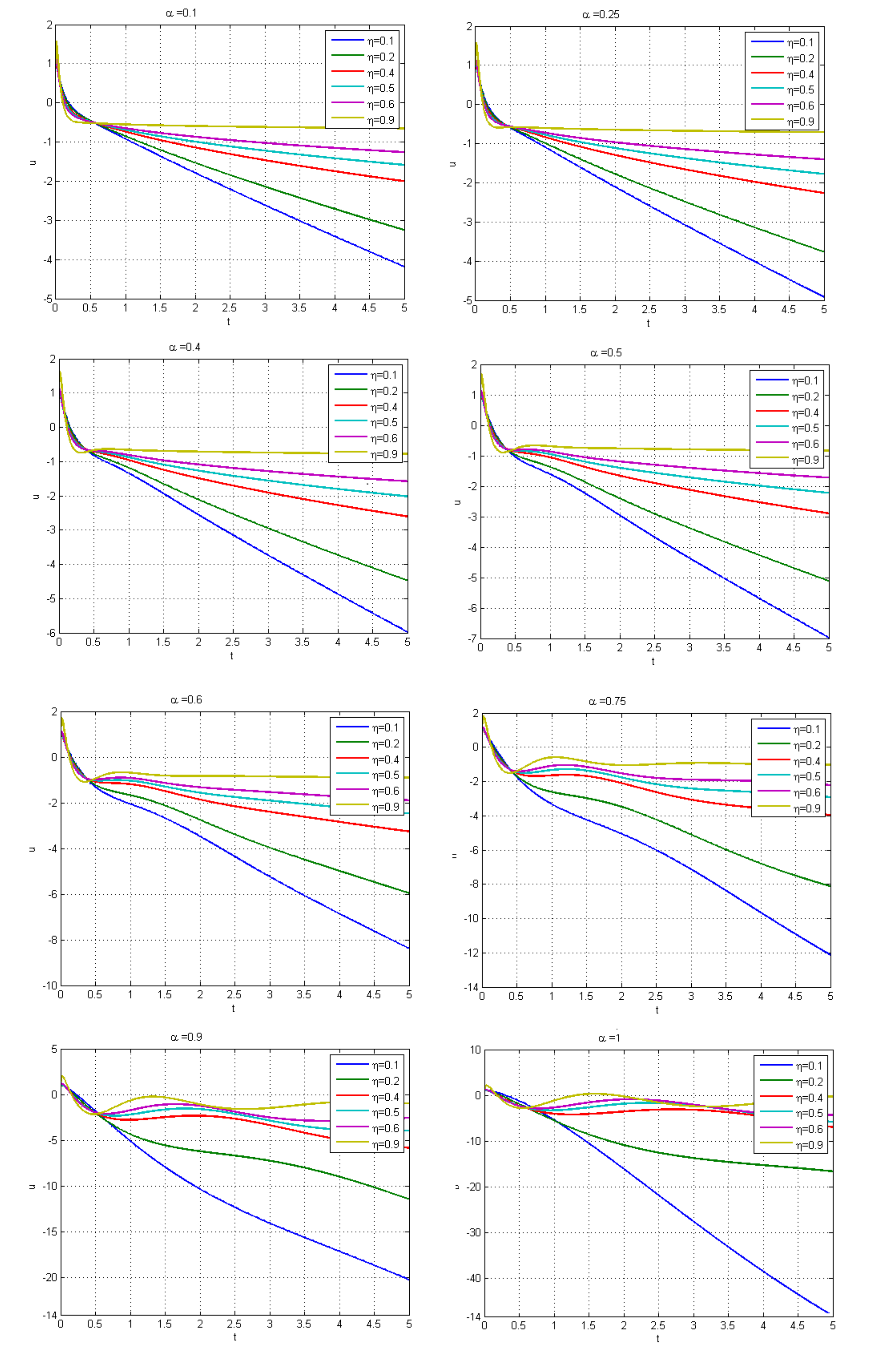}
\end{center}
\caption{Numerical solution of (\ref{hybrid_system}), plotted for different 
values of $\alpha \in (0,1]$ and $\eta$, with fixed $\lambda = 10$, initial 
discretization step $h = 0.001$, computed for $t \in [0,5]$.}
\end{figure} \label{more_alpha_fixed_lambda_h_001}

\smallskip

\section*{Acknowledgments}
This work was carried out within the framework of the Italian GNFM-INdAM.
%




\bigskip \smallskip

\it

\noindent

\noindent
$^1$ Department of Mathematics and Physics  \\
Roma Tre University \\
1, Largo S. Leonardo Murialdo \\
00146 Rome, ITALY  \\ 
e-mail: concezzi@mat.uniroma3.it \\
e-mail: spigler@mat.uniroma3.it 
 \\ 


\noindent
$^2$ Department of Statistical Sciences  \\
``Sapienza'' University of Rome  \\ 
5, P.le Aldo Moro \\
00185 Rome, ITALY \\ 
e-mail: roberto.garra@sbai.uniroma1.it  \\



\hfill Received: MONTH DD, YYYY
%


%
\end{document}

%% file: fcaa_style.tex
\usepackage{graphicx}
\usepackage{epsfig}
\usepackage{amsthm}
\usepackage{amsmath}
\usepackage{latexsym}
\usepackage{amsfonts}
\usepackage{amssymb}

\newtheoremstyle{theorem}
  {15pt}          
  {15pt}  
  {\sl}  
  {\parindent}
  {\sc}  
  {. }    
  { }    
  {}     
\theoremstyle{theorem}
\newtheorem{lemma}{Lemma}[section]

\newtheoremstyle{defi}
  {15pt}          
  {15pt}  
  {\rm}  
  {\parindent}     
  {\sc}  
  {. }    
  { }    
  {}     
\theoremstyle{defi}

 \def\proofname{\indent\rm P\ r\ o\ o\ f.}
 
 \def\qed{\hfill$\Box$} 
